\documentclass[10pt]{amsart}

\usepackage{setspace}
\usepackage{moreverb}

\newcommand{\ignore}[1]{}
\newenvironment{Macaulay2}{
\begin{spacing}{0.8}
\begin{quote}
\smallskip
\smallskip } { \smallskip
\end{quote}
\end{spacing}
\medskip
}

\newtheorem{theorem}{Theorem}[section]
\newtheorem{lemma}[theorem]{Lemma}

\newtheorem{algorithm}[theorem]{Algorithm}
\theoremstyle{definition}

\newtheorem{example}[theorem]{Example}

\theoremstyle{remark}
\newtheorem{remark}[theorem]{Remark}

\numberwithin{equation}{section}

\input xy
\xyoption{all}

\newcommand{\p}{\partial}
\newcommand{\grad}{\nabla}

\newcommand{\bC}{\mathbb{C}}

\newcommand{\bZ}{{\mathbb{Z}}}
\newcommand{\bQ}{{\mathbb{Q}}}

\newcommand{\rad}{\mbox{\rm{rad}} \,}
\newcommand{\Ann}{\mbox{\rm{Ann}} \,}

\newcommand{\dm}{\mathrm {dim }}

\newcommand{\Supp}{\mbox{\rm{Supp}} }

\newcommand{\spec}{\mathrm {Spec\, }}

\newcommand{\ckr}{\mathrm {Coker\, }}
\newcommand{\im}{\mathrm {Im\, }}

\newcommand{\Ext}{\mbox{\rm{Ext}} }

\newcommand{\fM}{{\mathfrak m}}

\newcommand{\fQ}{\mathfrak{q}}

\newcommand{\cA}{{\mathcal A}}

\newcommand{\cD}{{\mathcal D}}

\newcommand{\la}{{\lambda}}

\newcommand{\lra}{{\longrightarrow}}

\newcommand{\pa}{{\partial}}
\newcommand{\ta}{{T^{\ast} X}}

\begin{document}

\title[Computing the support of local cohomology modules]{Computing the support of local cohomology modules}

\author[J. \`Alvarez Montaner]{Josep \`Alvarez Montaner}
\thanks{Research of the first author partially supported by a
Fulbright grant and the Secretar\'{\i}a de Estado de Educaci\'on y Universidades of Spain and the European Social Funding}
\address{Departament de Matem\`atica Aplicada I\\
Universitat Polit\`ecnica de Catalunya\\ Avinguda Diagonal 647, Barcelona 08028, SPAIN} \email{Josep.Alvarez@upc.es}

\author[Anton Leykin]{Anton Leykin}
%\thanks{}
\address{Department
of Mathematics, Statistics, and Computer Science, University of Illinois at Chicago, 851 South Morgan (M/C 249), Chicago, IL
60607-7045, USA} \email{leykin@math.uic.edu}



\begin{abstract}
For a polynomial ring $R=k[x_1,...,x_n]$, we present a method to
compute the characteristic cycle of the localization $R_f$ for any
nonzero polynomial $f\in R$ that avoids a direct computation of
$R_f$ as a $D$-module. Based on this approach, we develop an
algorithm for computing the characteristic cycle of the local
cohomology modules $H^r_I(R)$ for any ideal $I\subseteq R$ using
the \v{C}ech complex. The algorithm, in particular, is useful for
answering questions regarding vanishing of local cohomology
modules and computing Lyubeznik numbers. These applications are
illustrated by examples of computations using our implementation
of the algorithm in Macaulay~2.
\end{abstract}

\maketitle

\section{Introduction}

Let $k$ be a field of characteristic zero and $R=k[x_1,\dots,x_n]$ the ring of polynomials in $n$ variables. For any ideal
$I\subseteq R$, the local cohomology modules $H_I^r(R)$ have a natural finitely generated module structure over the Weyl
algebra $A_n$. Recently, there has been an effort made towards effective computation of these modules by using the theory of
Gr\"obner bases over rings of differential operators. Algorithms given by U.~Walther \cite{Wa99} and T.~Oaku and N.~Takayama
\cite{OT01} provide a utility for such computation and are both implemented in the package {\tt D-modules} \cite{LST} for {\tt
Macaulay 2} \cite{GS}.

\vskip 2mm

Walther's algorithm is based on the construction of the $\check {\rm C}$ech complex in the category of $A_n$-modules. So it is
necessary to give a description of the localization $R_f$ at a polynomial $f\in R$. An algorithm to compute these modules was
given by T.~Oaku in \cite{Oa97}. The main ingredient of the algorithm is the computation of the Bernstein-Sato polynomial of
$f$ which turns out to be a major bottleneck due to its complexity.

\vskip 2mm

To give a presentation of $R_f$ or $H_I^r(R)$ as $A_n$-modules is out of the scope of this work. Our aim is to provide an algorithm to
compute an invariant that can be associated to a finitely generated $A_n$-module, the characteristic cycle. This invariant
gives a description of the support of the $A_n$-module as an $R$-module, so it is a useful tool to prove the vanishing of local
cohomology modules. Moreover, the characteristic cycles of local cohomology modules also give us some extra information since their
multiplicities are a set of numerical invariants of the quotient ring $R/I$ (see \cite{Al02}). Among these invariants we may
find Lyubeznik numbers that were introduced in \cite{Ly93}.

\vskip 2mm

We will present an algorithm to compute the characteristic cycle of
any local cohomology module. It comes naturally from the structure
of the \v{C}ech complex and the additivity of the characteristic
cycle with respect to short exact sequences. The requirement is that
we have to compute first the characteristic cycle of the
localizations appearing in the \v{C}ech complex. To do so, we
present a method based on a geometric formula given by V.~Ginsburg
in \cite{Gi86} and reinterpreted by J.~Brian\c{c}on, P.~Maisonobe
and M.~Merle in \cite{BMM94}. The advantage of this approach is that
we will not have to compute the Bernstein-Sato polynomial of $f$ and
we will be operating in a commutative graded ring in $2n$ variables
instead of operating in the Weyl algebra $A_n$. The algorithm we
will present is an elaboration of \cite[Thm. 3.4.2]{BMM94}, which --
we have to point out -- is stated in the complex analytic context.
For our computational purposes we are interested in the algebraic
context since we will need (absolute) primary decomposition so we
have to make sure that, at least for the examples we will develop,
this result is also true in the algebraic counterpart. The complex
algebraic case may already be found in \cite{Gi86} but the approach
we will use in this work is through flat base change. It allows us
to work over any field of characteristic zero for a large sample of
examples since localization modules and local cohomology modules
have a good behavior with respect to this operation. Since absolute
primary decomposition is not implemented in {\tt Macaulay 2}, we
compute over the field of rational numbers. This replacement,
however, is not an issue for the examples that we present.

\vskip 2mm

The scripts of the source codes we will use in this work as well as the output in full detail of the examples are available at
the web page http://www2.math.uic.edu/$\sim$leykin/CC.
In the future, we would explore the possibility of using numerical
primary decomposition -- technique for varieties over $\bC$ which is
being developed by the second author.

\section{Basics on the theory of $\cD$-modules}

Let $X=\bC^n$ be the complex analytic space with coordinate system
$x_1,\dots,x_n$. Given the convergent series ring
$R=\bC\{x_1,\dots,x_n\}$ consider the associated ring of
differential operators $D_n:=R\langle \pa_1,\dots,\pa_n\rangle$,
i.e. the ring extension generated by the partial derivatives
$\pa_i=\frac{\pa}{\pa x_i}$, with the relations given by
$\pa_i\pa_j=\pa_j\pa_i$ and $\pa_i r - r \pa_i=\frac{\pa}{\pa x_i}$,
where $ r\in R$. For any unexplained terminology concerning the
theory of rings of differential operators we shall use \cite{Bj79},
\cite{Co95}.

\vskip 2mm

The ring $D_n$ has a natural increasing filtration given by the total order; the corresponding associated graded ring $gr(D_n)$
is isomorphic to the polynomial ring $R[a_1,\dots,a_n]$. A finitely generated $D_n$-module $M$ has an increasing sequence of
finitely generated $R$-submodules such that the associated graded module $gr(M)$ is a finitely generated $gr(D_n)$-module. The
{ characteristic ideal} of $M$ is the ideal in $gr(D_n)=R[a_1,\dots,a_n]$ given by the radical ideal $J(M):= \rad (
\Ann_{gr(D_n)} (gr(M)))$. The ideal $J(M)$ is independent of the good filtration on $M$. The { characteristic variety} of $M$
is the closed algebraic set given by: $$C(M):= V(J(M))\subseteq \spec (gr(D_n))=\spec (R[a_1,\dots,a_n]).$$ The characteristic
variety describes the support of a finitely generated $D_n$-module as $R$-module. Let $\pi:
\spec(R[a_1,\dots,a_n])\longrightarrow \spec(R)$ be the map defined by $\pi(x,a)= x$. Then $\Supp_R(M)=\pi(C(M)).$

\vskip 2mm

We single out the important class of regular holonomic $D_n$-modules. Namely, a finitely generated $D_n$-module $M$ is
holonomic if $M=0$ or $\dim C(M)=n$. It is regular if there exists a good filtration on $M$ such that $\Ann_{gr(D_n)} (gr(M))$
is a radical ideal (\cite{BK}, see also \cite[\S 3]{Gi86}, \cite{Co95}).

\vskip 2mm

The { characteristic cycle} of $M$ is defined as:$$CC(M)= \sum m_i \hskip 2mm \Lambda_i$$ where the sum is taken over all the
irreducible components $\Lambda_i=V(\fQ_i)$ of the characteristic variety $C(M)$, where $\fQ_i \in \spec (gr(D_n))$ and $m_i$
is the multiplicity of $gr(M)$ at a generic point along each component $\Lambda_i$. These multiplicities can be computed via
Hilbert functions (see \cite{Ca84},\cite{Le}). Notice that the contraction of $\fQ_i$ to $R$ is a prime ideal so the variety
$\pi(\Lambda_i)$ is irreducible. These components can be described in terms of conormal bundles to $X_i:=\pi(\Lambda_i)$ in
$X$, i.e.  $$CC(M)=\sum m_i \hskip 2mm T_{X_i}^*X.$$ In particular, the support of $M$  is $\Supp_{R}(M)= \bigcup X_i $. For
details we refer, among others, to \cite[\S10]{Ph79}, \cite[\S7.5]{Ki88}.

\subsection{Characteristic cycle of a localization}

Let $M$ be a regular holonomic $D_n$-module. Then the localization
$M_f$ at a polynomial $f\in R$ is a regular holonomic $D_n$-module
as well. A geometric formula that provides the characteristic
cycle of $M_f$ in terms of the characteristic cycle of $M$ is
given by V.~Ginsburg in \cite{Gi86} and became known to us through
the interpretation of J.~Brian\c{c}on, P.~Maisonobe and M.~Merle
in~\cite{BMM94}.

\vskip 2mm

First we will recall how to compute the conormal bundle relative to $f$. Let $Y^\circ$ be the smooth part of a subvariety
$Y\subseteq X$ where $f|_{Y}$ is a submersion. Set:
$$W=\{(x,a)\in T^{\ast}X \hskip 2mm | \hskip 2mm x \in
Y^\circ \hskip 2mm {\rm and} \hskip 2mm a \hskip 2mm {\rm annihilates} \hskip 2mm T_x(f|_{Y})^{-1}(f(x))\}.$$ The conormal
bundle relative to $f$, denoted by $T_{f|_{Y}}^{\ast}$, is then the closure of $W$ in $\ta|_{Y}$.

\begin{theorem}\label{propBMM}{\rm (\cite[Thm. 3.4.2]{BMM94})}
Let $M$ be a regular holonomic $D_n$-module with characteristic cycle $CC(M)=\sum_i m_i \hskip 2mm T_{X_i}^*X$ and let $f\in R$
be a polynomial. Then $$CC(M_f)=\sum_{f(X_i)\neq 0} m_i(\Gamma_i+T_{X_i}^*X)$$ with $\Gamma_i=\sum_j m_{ij} \Gamma_{ij}$, where
$\Gamma_{ij}$ are the irreducible components of the divisor defined by $f$ in $T^*_{f|_{X_i}}$ and $m_{ij}$ are the
corresponding multiplicities.
\end{theorem}

% of the defining ideal of $\pi(\Gamma_{ij})$ along $\Gamma_{ij}$

\begin{remark}
Assume for simplicity that $M$ is a regular holonomic $D_n$-module such that $CC(M)= T_{Y}^*X$ and let $f\in R$ be a polynomial
such that $f(Y)\neq 0$. By the formula above we have $CC(M_f)=T_{Y}^*X + \Gamma$. It is worthwhile to point out that the
reduced variety associated to $\Gamma$ is the characteristic variety of the local cohomology module $H_{(f)}^1(M)$.

\end{remark}

\begin{example}
Set $R=\bC\{x,y,z\}$, $M=H^1_{(x)}(R)$, $f=x$ and $g=y$.

\vskip 2mm

%\begin{itemize}
%\item
 We have $CC(R)=T^*_{X}X$. Then
$T_{f|_{X}}^{\ast}= \{(x,y,z,a,b,c)\in T^{\ast}X \hskip 2mm | \hskip 2mm  b=0,c=0\}$
and the divisor defined by $f$ in $T^*_{f|_{X}}$ is
$\Gamma = \{(x,y,z,a,b,c)\in T^{\ast}X \hskip 2mm | \hskip 2mm  b=0,c=0,x=0\}=T^*_{\{x=0\}}X$. Thus $$CC(R_x)=T^*_{X}X+T^*_{\{x=0\}}X$$

\vskip 2mm
%\item
 We have $CC(M)=T^*_{\{x=0\}}X$. Then
$T_{g|_{\{x=0\}}}^{\ast}= \{(x,y,z,a,b,c)\in T^{\ast}X \hskip 2mm | \hskip 2mm  c=0,x=0\}$
and the divisor defined by $g$ in $T^*_{g|_{\{x=0\}}}$ is
$\Gamma = \{(x,y,z,a,b,c)\in T^{\ast}X \hskip 2mm | \hskip 2mm  c=0,x=0,y=0\}=T^*_{\{x=y=0\}}X$. Thus $$CC(M_y)=T^*_{\{x=0\}}X+T^*_{\{x=y=0\}}X$$
%\end{itemize}

\end{example}

The multiplicities $m_{ij}$ appearing in the formula are the
multiplicities of a generic point $x$ along each component
$\Gamma_{ij}$ of $\Gamma_i$ and can be computed via Hilbert
functions as in \cite{Le}.

\begin{lemma}\label{mult}
Let $e(\Gamma,x)$ denote the multiplicity of the variety $\Gamma
\subseteq T^*X$ defined by the ideal $I\subseteq R[a_1,\dots,a_n]$
at a point $x$. Then, the multiplicity $m$ of a generic point $x$
along $\Gamma$ is $$m=e(\Gamma,x)/e(\Gamma^{red},x),$$ where
$\Gamma^{red}$ is the variety defined by $\rad(I)$.

\end{lemma}

\begin{proof}

A reformulation of \cite[Prop. 3.11]{HIO} for the particular case of $x$ being a generic point gives us the desired result,
i.e. $e(\Gamma,x)= e(\Gamma^{red},x)\cdot m$.
\end{proof}

%The multiplicities of the characteristic cycle of a $D_n$-module $M$ can be computed via Hilbert functions on the ideal
%$\Ann_{gr(D_n)} (gr(M))$, i.e. it is not necessary to use the characteristic ideal $J(M)$ (see \cite{Le}).

\subsection{Algebraic $\cD$-modules}

Let $X=\bC^n$ be the complex affine space with coordinate system
$x_1,\dots,x_n$. Given the polynomial ring $R=\bC[x_1,\dots,x_n]$
consider the associated ring of differential operators
$A_n:=R\langle \pa_1,\dots,\pa_n\rangle$, i.e. the Weyl algebra. The
theories of algebraic $\cD$-modules and analytic $\cD$-modules are
very closely related. If one mimics the constructions given for the
ring $D_n$, one can check that the results we have considered
before, conveniently reformulated, remain true for $A_n$. In
particular we may construct an algebraic characteristic cycle as a
counterpart to the analytic characteristic cycle described before.
Our aim is to explain how both cycles are related.

\vskip 2mm

Set $\bC\{x\}:=\bC\{x_1,\dots,x_n\}$ and $\bC[x]:=\bC[x_1,\dots,x_n]$. Let $M$ be a regular holonomic $A_n$-module. The
$D_n$-module $M^{an}:=\bC\{x\}\otimes_{\bC[x]}M$ is also regular holonomic. For a good filtration $\{M_i\}_{i\geq 0}$ on $M$
the filtration $\{M_i^{an}:= \bC\{x\}\otimes_{\bC[x]}M_i\}_{i\geq 0}$ is also good due to the fact that $\bC\{x\}$ is flat over
${\bC[x]}$. Therefore $gr(M^{an})\simeq \bC\{x\}\otimes_{\bC[x]} gr(M)$ so the characteristic variety of $M^{an}$ is the
extension of the characteristic variety of $M$, i.e. $C(M^{an})= C(M)^{an}$. However, we should notice that the components of
the characteristic variety may differ depending on the ring we are considering. In particular we may have algebraically
irreducible components that are analytically reducible.

\vskip 2mm

The regular holonomic $A_n$-modules we will consider in this work, i.e.
the polynomial ring $R=\bC[x]$, the localization $R_f$ for a
polynomial $f\in R$, and the local cohomology modules $H_I^r(R)$,
all have a good behavior with respect to flat base change. We state
that, roughly speaking, the formulas of the algebraic and analytical
characteristic variety of these modules are the same but the
components and  multiplicities of the corresponding characteristic
cycle may differ.

\begin{remark}
The results of this section can be stated in general for $X$ being any smooth algebraic variety over $\bC$. It is worth to
point out that $M\rightarrow M^{an}$ gives an equivalence between the category of regular holonomic $\cD_X$-modules and the
category of regular holonomic $\cD_X^{an}$-modules when $X$ is projective (see \cite[\S 3]{Gi86}).

\end{remark}

\section{Algorithmic approach to Brian\c{c}on-Maisonobe-Merle's Theorem}

From now on we will assume that $R=\bC[x_1,\dots,x_n]$ is the
polynomial ring so we will be working in the algebraic context. Let
$M$ be a regular holonomic $A_n$-module with algebraic
characteristic cycle $CC(M)=\sum m_i \hskip 2mm T_{X_i}^*X$ and let
$f\in R$ be a polynomial. Our aim is to compute the characteristic
cycle of the localization $M_f$ operating in the commutative graded
ring $gr(A_n)=R[a_1,\dots,a_n]$.  We are going to provide two
algorithms that are an elaboration of Theorem \ref{propBMM}. The
first one computes the part $\Gamma_i$ of the formula in Theorem
\ref{propBMM} corresponding to each irreducible component
$T_{X_i}^*X$ in the characteristic cycle of $M$. The second one
computes the components and the corresponding multiplicities of the
varieties $\Gamma_i$.

\vskip 2mm

Theorem \ref{propBMM} is a geometric reformulation of a result given by V.~Ginsburg \cite[Thm. 3.3]{Gi86}. Even though it is
stated in the analytic context me may find in Ginsburg's paper the algebraic counterpart to the same result, see \cite[Thm.
3.2]{Gi86}. We may interpret it as in Section 2.2. through flat base change.

\vskip 2mm

\begin{algorithm} \label{alg1}{ (Divisor defined by $f$ in $T^*_{f|_{Y}}$, the conormal relative to $f$)}
\vskip 2mm

{\rm

\noindent {\sc Input:} Generators $g_1,...,g_d$ of an ideal $I\subset R$ defining the algebraic variety $Y=V(I)\subseteq X$ and
a polynomial $f\in R$.

\noindent {\sc Output:} Divisor defined by $f$ in the conormal $T^*_{f|_{Y}}$ relative to $f$.

\vskip 2mm

{\bf Compute the smooth part $Y^\circ$ of $Y$ where $f|_{Y}$ is a submersion:}

\begin{itemize}

\item [\textbf{(0a)}] Compute $\grad f=(\frac{\p f}{\p
x_1},...,\frac{\p f}{\p x_n})$

\item [\textbf{(0b)}] Compute  $Y^{\circ} = Y\setminus V(I^\circ)$, where $I^\circ\subset R$ is the defining ideal of $ \{x \in Y \ | \ \grad f(x)=0 \}$  and the singular locus of $Y$.

%\item [\textbf{(0b)}] Compute the ideal $I^\circ\subset R$ such
%that $Y^{\circ} = \{x \in Y \ | \ \grad f(x)\notin T_x Y \}$ is described as $Y^{\circ} = Y\setminus V(I^\circ)$.

%\marginpar{ In principle there could be a situation when $X^\circ$
%is empty! Example: let $X$ be a line colinear with the gradient.}

\end{itemize}

{\bf Compute the conormal relative to $f$}

\begin{itemize}
\item[\textbf{(1a)}] Compute $K = \ker \phi$, where the $\phi:
R^n\to R^{d+1}/I$ sends
$$
s\mapsto (\nabla f, \nabla g_1, ..., \nabla g_d)\cdot s \in R^{d+1}/I.
$$

\item[\textbf{(1b)}] Let $J\subset gr(A_n) = R[a_1,...,a_n]$ be
the ideal generated by $\{ (a_1,...,a_n)\cdot b \ |\ b \in K \}.$

\item[\textbf{(1c)}]  Compute $J_{sat}=J:(gr(A_n)I^\circ)^\infty$. (Note: $I(T^*_{f|_{Y}}) =
\sqrt{J_{sat}}$.)

\end{itemize}

{\bf Compute the divisor defined by $f$ in $T^*_{f|_{Y}}$ }

\begin{itemize}

\item[\textbf{(2a)}] Compute $K_f = \ker \phi_f$, where the map
$\phi_f: R^n\to R^{d+1}/(I+(f))$ sends
$$
s\mapsto (\nabla f, \nabla g_1, ..., \nabla g_d)\cdot s \in R^{d+1}/(I+(f)).
$$

\item[\textbf{(2b)}] Let $J_f \subset gr(A_n) = R[a_1,...,a_n]$ be
the ideal generated by $\{ (a_1,...,a_n)\cdot b \ |\ b \in K_f \}$.

\item[\textbf{(2c)}] $C = J_{sat} + (f) + J_f\subset gr(A_n)$.
\end{itemize}

\noindent {\sc Return:} The ideal $C$ that defines the divisor $f$ in $T^*_{f|_{Y}}$ }
\end{algorithm}

\begin{proof}(Correctness of the algorithm)
The steps (0a), (0b) follow from the definition of $f|_{Y}$ being a submersion. The relative conormal $T^*_{f|_{Y}}$ is the
closure of $$W=\{(x,a)\in T^*X \ | \ x\in Y^\circ , \forall s\in K,\ a(s)=0 \}.$$ For every
point $ x\in Y^\circ$, the tangent space $T_{x}Y^\circ$ is a specialization of $V(K)$, where $K$ is computed in step (1a). A
defining ideal of $W$ is produced in (1b) and, finally, taking the closure amounts to the saturation in (1c). In order to
restrict to $f=0$, it is not enough to compute $J_{sat} + (f)$. However, step (2a) and (2b) that follow closely the idea of
(1a) and (1b) provide the necessary correction term in (2c).

\vskip 2mm

Recall that the analytic extension of the ideal $C$ we obtain with the algorithm is what we would obtain applying Theorem
\ref{propBMM} in order to compute the analytic characteristic cycle of the localization module (see Section 2.2). In our case,
the ideal $C$ will give us the components of the algebraic characteristic variety.

\end{proof}

\begin{algorithm} \label{alg2} { (Components and multiplicities of the characteristic cycle)}

\vskip 2mm

{\rm \noindent {\sc Input:} The characteristic cycle $CC(M)=\sum m_i \hskip 2mm T_{X_i}^*X$ of a regular holonomic
$A_n$-module $M$ and a polynomial $f\in R$.

\noindent {\sc Output:} The characteristic cycle $CC(M_f)=\sum_{f(X_i)\neq 0} m_i(\Gamma_i+T_{X_i}^*X)$.

\vskip 2mm

For every component $Y=X_i$ we have to compute the ideal $C_i$ corresponding to the divisor defined by $f$ in $T^*_{f|_{Y}}$
using Algorithm \ref{alg1}. Then:

\vskip 2mm

{\bf Compute the components of $C_i$}

\begin{itemize}

\item[\textbf{(1a)}] Compute the associated primes $C_{ij}$ of $C_i$.

\item[\textbf{(1b)}] Compute $I_{ij}= C_{ij} \cap R$ (if you need to know the defining
ideal of $X_{ij}= \pi(\Gamma_{ij})$ in Theorem \ref{propBMM}).

\end{itemize}

{\bf Compute the multiplicities}

\begin{itemize}

\item[\textbf{(2)}] Compute the multiplicity $m_{ij}$ in Theorem \ref{propBMM} as the multiplicity of a generic point $x$ along each
component $C_{ij}$ of $C_i$ as in Lemma \ref{mult}  via Hilbert functions.

%\marginpar{\small Have no good way of computing multiplicities in the general case.}

\end{itemize}

\noindent {\sc Return:} The components of $CC(M_f)$ and their corresponding multiplicities.}

\end{algorithm}

\begin{proof} The correctness of the algorithm is straightforward and follows from
Lemma \ref{mult}.
\end{proof}

%\subsection{Implementation of the algorithm}

The algorithm we propose requires the computation of the
associated primes of an ideal; primary decomposition is also
needed in the implementation if we want to avoid choosing generic
points when computing the multiplicities (see Lemma \ref{mult}).
Therefore, we have to restrict ourselves to computations in the
polynomial ring $R=\bQ[x_1,\dots,x_n]$ as we implemented the
algorithm in the computer system {\tt Macaulay~2}. What we are
going to construct is the characteristic cycle of a regular
holonomic $A_n$-module where now $A_n$ stands for the Weyl algebra
with rational coefficients. By flat base change we can extend the
ideal $C$ we obtain with Algorithm \ref{alg1} to any ring of
polynomials over a field of characteristic zero or to the
convergent series ring over $\bC$. As we stated in Section 2.2,
the primary components may differ depending on the ring we are
considering.

\vskip 2mm

In order to construct the algebraic characteristic cycle over $\bQ$
we would need to find the absolute primary decomposition of the
ideal $C$ we obtain with Algorithm \ref{alg1}. Even though the {\tt
Macaulay 2} command for primary decomposition is not implemented
over the algebraic closure of $\bQ$, it suffices for the examples we
treat in the next section.

\vskip 2mm

Another fine point in the implementation is the treatment of
embedded components of the ideal $C$ outputted in Algorithm
\ref{alg1}. The ideal $C$ contains the complete information about
maximal components of the divisor $f$ on $T^*_{f|_{Y}}$, in
particular, we can compute their multiplicities. However, the
primary ideals in the decomposition of $C$ that correspond to an
embedded component may not lead to the correct multiplicity due to
the global nature of our computations. In order to obtain this
multiplicity we restrict the divisor to the embedded component,
which amounts to rerunning Algorithm \ref{alg1} `modulo' its
defining ideal. The top-level routine in our implementation
processes components recursively ``descending'' to, i.e., localizing
at, the embedded components when needed.

\section{Characteristic cycle and \v{C}ech complex}

Let $I=(f_1,\dots,f_s)\subseteq R=\bC[x_1,\dots,x_n]$  be an ideal and $M$ be a holonomic $A_n$-module. In this section we are going to compute the
characteristic cycle of the local cohomology modules $H^r_I(M)$ using the \v{C}ech complex
$$
{\check{C}}^{\bullet}(f_1,\dots,f_s;M) : \hskip 5mm 0 \longrightarrow M \stackrel{d_0}\longrightarrow \bigoplus_{i=1}^s M_{f_i} \stackrel{d_1}\longrightarrow
\cdots \longrightarrow M_{f_1\cdots f_s}\longrightarrow 0.
$$

For simplicity we will assume from now on that $M$ is
indecomposable. Otherwise, if $M=M_1 \oplus M_2$, then $$
{\check{C}}^{\bullet}(f_1,\dots,f_s;M) =
{\check{C}}^{\bullet}(f_1,\dots,f_s;M_1) \oplus
{\check{C}}^{\bullet}(f_1,\dots,f_s;M_2)$$ and $ H^r_I(M)=
H^r_I(M_1)\oplus H^r_I(M_2)$ for all $r$, so we can compute the
characteristic cycle of both local cohomology modules separately.
Sometimes we will denote the localization modules appearing in the
\v{C}ech complex $M_{f_{\alpha}}$, where
$f_{\alpha}=\prod_{\alpha_i=1} f_i$ for all $\alpha\in \{0,1\}^s$.
We will also denote $|\alpha|= \alpha_1 +\cdots+\alpha_s$ and
$\varepsilon_1,\dots, \varepsilon_s$ will be the natural basis of
$\bZ^s$.

\vskip 2mm

>From the characteristic cycle of the localization modules in the
complex we develop an algorithm to extract the precise information
needed to describe the characteristic cycles of the local cohomology
modules. The algorithm comes naturally from the structure of the
\v{C}ech complex and the additivity of the characteristic cycle with
respect to short exact sequences. However the following assumption
will be required:

\vskip 2mm

$(\dagger)$ For all $\alpha\in \{0,1\}^s$ such that $\alpha_i=0$,
the localization map $M_{f_{\alpha}} \lra M_{f_{\alpha +
\varepsilon_i}}$ is either a natural inclusion, i.e.,
$M_{f_{\alpha}}$ is saturated with respect to $f_{\varepsilon_i}$,
or $M_{f_{\alpha + \varepsilon_i}}=0$.

\vskip 2mm For unexplained terminology on the theory of complexes we
refer to \cite{Ro}. To shed some light on the process we first
present the case of $I$ being generated by one and two elements.

\subsection{The case $s=1$} We have the short exact sequence
$$
0 \longrightarrow H^0_I(M) \longrightarrow M \stackrel{d_0}\longrightarrow M_{f_1}
\longrightarrow H^1_I(M)\lra 0
$$
Under the assumption $(\dagger)$ either $CC(H^1_I(M))=
CC(M_{f_1})-CC(M)$ if $M$ is $f_1$-saturated or $CC(H^0_I(M))=
CC(M)$ if $M_{f_1}= 0$. The algorithm we propose to compute the
characteristic cycle works in both cases and boils down to the
following step:

\vskip 2mm

{\it Prune} the characteristic cycles of $M_{f_1}$ and $M$, where
`prune' means remove the components (counting multiplicities) that
both modules have in common.

\vskip 2mm

The characteristic cycle of $H^0_I(R)$ (resp. $H^1_I(R)$) is the
formal sum of components of $M$ (resp. $M_{f_1}$) that survived this
process.

\subsection{ The case $s=2$} The vertical sequences of the following diagram are exact
$${\xymatrix { {\check{C}}^{\bullet}(f_1;M): & 0 \ar[r] & M\ar[r] & M_{f_1}\ar[r] & 0 &
\\ {\check{C}}^{\bullet}(f_1,f_2;M): & 0 \ar[r] & M \ar[u] \ar[r]  & M_{f_1} \oplus M_{f_2} \ar[u] \ar[r] & M_{f_1f_2} \ar[r]& 0
\\ {\check{C}}^{\bullet}(f_1;M_{f_2})[-1]: & & 0 \ar[r] & M_{f_2} \ar[u] \ar[r]&M_{f_1f_2} \ar[u] \ar[r]& 0 }}
$$ so we have an exact sequence of \v{C}ech complexes $$ (i) \hskip 5mm 0 \longrightarrow {\check{C}}^{\bullet}(f_1;M_{f_2})[-1]
\longrightarrow {\check{C}}^{\bullet}(f_1,f_2;M) \longrightarrow {\check{C}}^{\bullet}(f_1;M)\lra 0$$ where $[-1]$ stands for the result of shifting the complex one place to the right. Analogously, i.e. switching the roles played by $f_1$ and $f_2$, we also have   $$ (ii) \hskip 5mm 0 \longrightarrow {\check{C}}^{\bullet}(f_2;M_{f_1})[-1]
\longrightarrow {\check{C}}^{\bullet}(f_1,f_2;M) \longrightarrow {\check{C}}^{\bullet}(f_2;M)\lra 0$$

\vskip 2mm

Notice that the vanishing of any localization module reduces the \v{C}ech complex to the case $s=1$. So we are going to consider the only case remaining under the assumption $(\dagger)$, i.e. $M$ is saturated with respect to $f_1f_2$.
Consider the  long exact sequence
of cohomology modules associated to $(i)$
$$
0  \lra H^{-1}_{(f_1)}(M_{f_2}) \longrightarrow H^0_I(M)
\longrightarrow H^0_{(f_1)}(M)\stackrel{\delta^0}\lra
H^0_{(f_1)}(M_{f_2}) \longrightarrow H^1_I(M) \longrightarrow
\cdots
$$
where $\delta^j$ are the connecting maps. Non-vanishing may occur
only in the sequence
$$
0 \longrightarrow H^1_I(M) \longrightarrow
H^1_{(f_1)}(M)\stackrel{\delta^1}\lra H^1_{(f_1)}(M_{f_2})
\longrightarrow H^2_I(M) \longrightarrow 0
$$
which breaks down into two short exact sequences with $C_1 = \ckr
\delta^1$:

\vskip 2mm

\hskip 5mm  $0\longrightarrow H^1_I(M) \longrightarrow H^1_{(f_1)}(M)\lra C_1 \longrightarrow 0 $ \hskip 2mm  and \hskip 2mm $ 0\longrightarrow C_1 \lra
H^1_{(f_1)}(M_{f_2}) \longrightarrow H^2_I(M) \longrightarrow 0$

\vskip 2mm

In order to get  $CC(H^1_I(M))$ and $CC(H^2_I(M))$ we only have to compute the characteristic cycle of $C_1$ since we already know that
$CC(H^1_{(f_1)}(M))=CC(M_{f_1})- CC(M)$  and $CC(H^1_{(f_1)}(M_{f_2}))=CC(M_{f_1f_2})- CC(M_{f_2})$.
\vskip 2mm

{\bf Claim:} $CC(C_1)= \sum m_i \hskip 1mm T_{X_i}^{\ast} X$ where the sum is taken over the components (counting multiplicities)
that $CC(H^1_{(f_1)}(M))$ and $CC(H^1_{(f_1)}(M_{f_2}))$ have in common.

\begin{proof}
Assume that there is a component $T_{X_i}^{\ast} X$  in
$CC(H^1_{(f_1)}(M))$ and $CC(H^1_{(f_1)}(M_{f_2}))$  not appearing
in $CC(C_1)$, i.e.  $T_{X_i}^{\ast} X$ is a component of
$CC(H^1_I(M))$ and $CC(H^2_I(M))$. This component shows up in the
computation of the characteristic cycle of the cohomology of the
subcomplex ${\check{C}}^{\bullet}(f_2;M_{f_1})[-1]$, since it is in
fact a component of $ CC(M_{f_1})$ and  $CC(M_{f_1f_2})$, but is not
a component of $CC(M)$ and $CC(M_{f_2})$. Consider the  long exact
sequence of cohomology modules associated to the complex $(ii)$:
$$  \cdots \lra H^0_I(M)
\longrightarrow H^0_{(f_2)}(M)\stackrel{\delta^0}\lra
H^0_{(f_2)}(M_{f_1}) \longrightarrow H^1_I(M) \longrightarrow
H^1_{(f_2)}(M) \longrightarrow \cdots
$$
It follows that the component $T_{X_i}^{\ast} X$ should belong to $CC(H^0_{(f_2)}(M_{f_1}))$ and  $CC(H^1_{(f_2)}(M_{f_1}))$ in order to fulfill the
hypothesis of being a component of $CC(H^1_I(M))$ and $CC(H^2_I(M))$. Thus we get a
contradiction since $H^0_{(f_2)}(M_{f_1})=0$ as $M_{f_1}$ has no
$f_2$-torsion.
\end{proof}

To summarize the computation of the characteristic cycle of
$H^r_I(M)$ using the sequence of \v{C}ech complexes $(i)$ we propose
the following algorithm:

\vskip 2mm

$(1)$ {\it Prune} the characteristic cycle of $M_{f_1}$ and $M$.

\hskip .6cm {\it Prune} the characteristic cycle of $M_{f_1f_2}$
and $M_{f_2}$.

\vskip 2mm

$(2)$ {\it Prune} the characteristic cycle of $M_{f_2}$ and $M$.

\hskip .6cm {\it Prune} the characteristic cycle of $M_{f_1f_2}$
and $M_{f_1}$.

\vskip 2mm

`Prune' means remove the components (counting multiplicities) such
that both characteristic cycles still have in common in that step of
the algorithm. The characteristic cycle of $H^0_I(M)$ (resp.
$H^1_I(M)$, $H^2_I(M)$) is the formal sum of comand  $CC(H^1_{(f_2)}(M_{f_1}))$ponents of $M$
(resp. $M_{f_1}$ and $M_{f_2}$, $M_{f_1f_2}$), i.e., components of
th characteristic cycles of ${\check{C}}^{0}(f_1,f_2;M)$ (resp.
${\check{C}}^{1}(f_1,f_2;M)$, ${\check{C}}^{2}(f_1,f_2;M)$) that
survived to the process.

\vskip 2mm

\begin{remark}
Naturally, one may permute the steps (1) and (2) in the above
algorithm.
\end{remark}

\subsection{The general case}
Let $I=(f_1,\dots,f_s)\subseteq R$ be an ideal and $M$ be a holonomic $A_n$-module satisfying

\vskip 2mm

$(\dagger)$ For all $\alpha\in \{0,1\}^s$ such that $\alpha_i=0$, the localization map $M_{f_{\alpha}} \lra M_{f_{\alpha + \varepsilon_i}}$ is a natural inclusion, i.e. $M_{f_{\alpha}}$ is saturated with respect to $f_{\varepsilon_i}$, or $M_{f_{\alpha + \varepsilon_i}}=0$.
\vskip 2mm

Our aim is to proceed inductively in order to compute the characteristic cycle of the local cohomology modules $H_I^r(M)$.
To this purpose it is useful to visualize the \v{C}ech complex ${\check{C}}^{\bullet}(f_1,\dots,f_s;M)$ as a $s$-hypercube where the edges are the localization maps that, with the corresponding sign, describe the differentials of the complex. The exact sequence of complexes
$$ 0 \longrightarrow {\check{C}}^{\bullet}(f_1,\dots,f_{s-1};M_{f_{s}})[-1] \longrightarrow
{\check{C}}^{\bullet}(f_1,\dots,f_{s};M) \longrightarrow {\check{C}}^{\bullet}(f_1,\dots,f_{s-1};M)\lra 0$$ can be easily identified in the $s$-hypercube. For the case $s=3$ we visualize the \v{C}ech complexes ${\check{C}}^{\bullet}(f_1,f_2,f_3;M)$,
${\check{C}}^{\bullet}(f_1,f_{2};M)$  and ${\check{C}}^{\bullet}(f_1,f_2;M_{f_{3}})[-1]$ as follows

{\tiny 
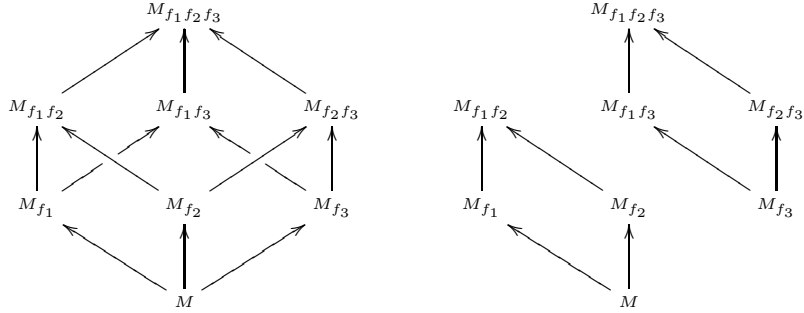
$${\xymatrix { &M_{f_1f_2f_3} &
\\ M_{f_1f_2} \ar[ur]  & M_{f_1f_3} \ar[u] & M_{f_2f_3}  \ar[ul]
\\ M_{f_1} \ar[ur]|\hole \ar[u]&M_{f_2} \ar[ul] \ar[ur]& M_{f_3} \ar[ul]|\hole \ar[u]
\\& M \ar[ul] \ar[u] \ar[ur]& }}   \hskip 1cm
{\xymatrix
{ &M_{f_1f_2f_3} &
\\ M_{f_1f_2}   & M_{f_1f_3} \ar[u]  & M_{f_2f_3} \ar[ul]
\\ M_{f_1}  \ar[u]& M_{f_2} \ar[ul] & M_{f_3} \ar[ul] \ar[u]
\\& M \ar[ul] \ar[u] & }}
$$}

Chasing the diagrams, one may check that the localization map with
respect to $f_s$, i.e., the edges $M_{f_\alpha}\lra M_{f_{\alpha
+\varepsilon_s}}$ in the $s$-hypercube, induces the connecting maps
$\delta^j$ in the long exact sequence of cohomology modules $$  0
\lra H^{-1}_{(f_1,...,f_{s-1})}(M_{f_s}) \longrightarrow H^0_I(M)
\longrightarrow H^0_{(f_1,...,f_{s-1})}(M)\stackrel{\delta^0}\lra
H^0_{(f_1,...,f_{s-1})}(M_{f_s}) \longrightarrow H^1_I(M)
\longrightarrow \cdots$$

We are not going to give a precise description of the connecting
maps, since we are only interested in the data given by the
characteristic cycle. The formula we obtain in Theorem \ref{T} for
the characteristic cycle of the local cohomology modules $H_I^r(M)$
is given just in terms of the components of the characteristic cycle
of the localizations $M_{f_{\alpha}}$. The precise information we
need to extract is given by the following algorithmic procedure.

\begin{algorithm} \label{alg3} { (Characteristic cycle and \v{C}ech complex )}

\vskip 2mm

{\rm \noindent {\sc Input:} Characteristic cycles
$CC(M_{f_{\alpha}})= \sum m_{\alpha,i} \hskip 1mm T_{X_i}^{\ast} X$
for $\alpha \in \{0,1\}^s$.

\noindent {\sc Output:} (Pruned) characteristic subcycles $\overline{CC}
(M_{f_{\alpha}}) \subseteq CC(M_{f_{\alpha}})$ for $\alpha \in
\{0,1\}^s$.

\vskip 2mm

{\bf Prune the extra components}

\vskip 2mm

For $j$ from $1$ to $s$, incrementing by $1$

\vskip 2mm

\begin{itemize}

\item[\textbf{(j)}] ${\it Prune}$ the localizations
$M_{f_{\alpha}} $ and $M_{f_{\alpha + \varepsilon_j}}$ for all
$\alpha \in \{0,1\}^s$ such that $\alpha_j=0$, where `prune' means
remove the components (counting multiplicities) such that both
modules still have in common after step $(j-1)$.

\end{itemize}

\vskip 2mm

\noindent {\sc Return:} The components and the corresponding
multiplicities of a characteristic subcycle of $CC(M_{f_{\alpha}})$
for $\alpha \in \{0,1\}^s$.

}

\end{algorithm}

For the case $s=3$ we visualize the steps of the algorithm  as follows \hskip -1cm{\tiny
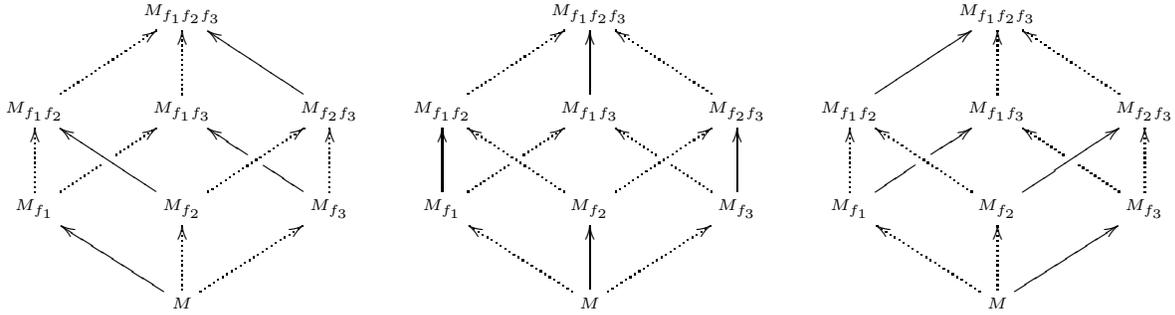
$${\xymatrix { &M_{f_1f_2f_3} &
\\ M_{f_1f_2} \ar@{.>}[ur]  & M_{f_1f_3} \ar@{.>}[u] & M_{f_2f_3}  \ar[ul]
\\ M_{f_1} \ar@{.>}[ur]|\hole \ar@{.>}[u]&M_{f_2} \ar[ul] \ar@{.>}[ur]& M_{f_3} \ar[ul]|\hole \ar@{.>}[u]
\\& M \ar[ul] \ar@{.>}[u] \ar@{.>}[ur]& }}   \hskip .51cm
{\xymatrix { &M_{f_1f_2f_3} &
\\ M_{f_1f_2} \ar@{.>}[ur]  & M_{f_1f_3} \ar[u] & M_{f_2f_3} \ar@{.>}[ul]
\\ M_{f_1} \ar@{.>}[ur]|\hole \ar[u]&M_{f_2} \ar@{.>}[ul] \ar@{.>}[ur]& M_{f_3} \ar@{.>}[ul]|\hole \ar[u]
\\& M \ar@{.>}[ul] \ar[u] \ar@{.>}[ur]& }}   \hskip .51cm
{\xymatrix { &M_{f_1f_2f_3} &
\\ M_{f_1f_2} \ar[ur]  & M_{f_1f_3} \ar@{.>}[u] & M_{f_2f_3} \ar@{.>}[ul]
\\ M_{f_1} \ar[ur]|\hole \ar@{.>}[u]&M_{f_2} \ar@{.>}[ul] \ar[ur]& M_{f_3} \ar@{.>}[ul]|\hole \ar@{.>}[u]
\\& M \ar@{.>}[ul] \ar@{.>}[u] \ar[ur]& }}
$$}
The solid arrows indicate the modules we must prune at each step.

\begin{remark}
{\rm  As in the case $s=2$,  the order we propose in the algorithm depends on the \v{C}ech subcomplexes we consider
when computing the characteristic cycle of the cohomology of the \v{C}ech complex. We can obtain equivalent algorithms permuting the generators of the ideal $I$.  It is also worth mentioning that the algorithm can be used in the algebraic context over any field of characteristic zero and in the analytic context.}

\end{remark}

\begin{theorem}\label{T}
Let $I=(f_1,\dots,f_s)\subseteq R$ be an ideal and $M$ be an indecomposable holonomic $A_n$-module satisfying
$(\dagger)$. Then
$$
CC(H_I^r(M))= \sum_{|\alpha|=r} \overline{CC}(M_{f_\alpha}),
$$
where the pruned characteristic subcycles are obtained with
Algorithm \ref{alg3}.

\end{theorem}

\begin{proof}
We proceed by induction on the number of generators $s$ of the ideal. The cases $s=1,2$ have been already done. For $s>2$ we
have an exact sequence of complexes
$$ 0 \longrightarrow {\check{C}}^{\bullet}(f_1,\dots,f_{s-1};M_{f_{s}})[-1] \longrightarrow
{\check{C}}^{\bullet}(f_1,\dots,f_{s};M) \longrightarrow {\check{C}}^{\bullet}(f_1,\dots,f_{s-1};M)\lra 0$$
Splitting the corresponding associated long exact sequence of cohomology modules $$  0  \lra H^{-1}_{(f_1,...,f_{s-1})}(M_{f_s}) \longrightarrow H^0_I(M)
\longrightarrow H^0_{(f_1,...,f_{s-1})}(M)\stackrel{\delta^0}\lra H^0_{(f_1,...,f_{s-1})}(M_{f_s}) \longrightarrow H^1_I(M) \longrightarrow \cdots$$ into short exact sequences we obtain

\vskip 2mm

\hskip 1cm $0\lra A_r \lra H^r_{I}(M)\lra B_r \lra 0$

\vskip 2mm

\hskip 1cm $0\lra B_r \lra H^r_{(f_1,\dots,f_{s-1})}(M)\lra C_{r} \lra 0$

\vskip 2mm

\hskip 1cm $0\lra C_{r} \lra
H^{r}_{(f_1,\dots,f_{s-1})}(M_{f_s})\lra A_{r+1} \lra 0$

\vskip 2mm

The characteristic cycle of $H^r_{(f_1,\dots,f_{s-1})}(M)$  (resp.
$H^{r}_{(f_1,\dots,f_{s-1})}(M_{f_s})$) is the formal sum of
components of $M_{f_{\alpha}}$ satisfying $\alpha_s=0$ and
$|\alpha|=r$ (resp. $\alpha_s=1$ and $|\alpha|=r+1$) that survived
to step $(s-1)$ of the algorithm. Thus, for every $r$, in order to
get $CC(H^r_{I}(M))$ we only have to compute $CC(C_r)$ and use
additivity of the characteristic cycle with respect to short exact
sequences.

\vskip 2mm

{\bf Claim:} $CC(C_r)= \sum m_i \hskip 1mm T_{X_i}^{\ast} X$ where the sum is taken over the components (counting
multiplicities) that $CC(H^r_{(f_1,\dots,f_{s-1})}(M))$ and $CC(H^{r}_{(f_1,\dots,f_{s-1})}(M_{f_s}))$ have in common.

\vskip 2mm

Assume that there is a component $T_{X_i}^{\ast} X$  in $CC(H^r_{(f_1,\dots,f_{s-1})}(M))$ and
$CC(H^{r}_{(f_1,\dots,f_{s-1})}(M_{f_s}))$  not appearing in $CC(C_r)$, i.e. it has not been pruned in step $(s)$. In particular, this component belongs to the characteristic cycle of some localization modules $M_{f_{\alpha}}$ and $M_{f_{\alpha+\varepsilon_s}}$ satisfying  $\alpha_s=0$ and $|\alpha|=r$.
Then, this component is not pruned in the computation of the characteristic cycle of the cohomology of a convenient proper \v{C}ech
subcomplex of $ {\check{C}}^{\bullet}(f_1,\dots,f_s;M)$ containing  $M_{f_{\alpha}}$ and $M_{f_{\alpha+\varepsilon_s}}$. Thus we get a contradiction.

\end{proof}

If $M$ is not indecomposable we only have to apply Theorem \ref{T}
to each component.

\begin{example}
Consider $R=\bC[x]$ and the holonomic $A_1$-module $M= R \oplus
H_{(x)}^1(R)$.
  We have:

\vskip 2mm

$CC(M)= T^*_{X}X + T^*_{\{x=0\}}X$

$CC(M_x)= T^*_{X}X + T^*_{\{x=0\}}X$.

\vskip 2mm

\noindent Applying the pruning algorithm to each component, that
satisfy $(\dagger)$, we get
$$
CC(H^0_{(x)}(M))= CC(H^1_{(x)}(M))=
T^*_{\{x=0\}}X.
$$

One may be tempted to apply the pruning algorithm to $M$, however,
it misleads us resulting in the seeming vanishing of the local
cohomology modules $H^r_{(x)}(M)$. Notice that $M$ is not saturated
with respect to $f=x$.

\end{example}

The question whether Theorem \ref{T} still holds for indecomposable $A_n$-modules not satisfying $(\dagger)$ is
open. For some examples we may give an affirmative answer.

\begin{example}
Set $R=\bC[x,y]$. The holonomic $A_2$-module $M=H_{(xy)}^1(R)$ is not saturated with respect to $f=y$. We have:

\vskip 2mm

$CC(M)= T^*_{\{x=0\}}X +T^*_{\{y=0\}}X + T^*_{\{x=y=0\}}X$

$CC(M_y)= T^*_{\{x=0\}}X + T^*_{\{x=y=0\}}X$.

\vskip 2mm

\noindent Pruning the components they have in common we get  $CC(H^0_{(y)}(M))= T^*_{\{y=0\}}X$. It agrees with the fact that $H^0_{(y)}(M)\cong H_{(y)}^1(R)$.

\end{example}

\begin{remark}

Let $I=(f_1,\dots,f_s)\subseteq R$ be an ideal and $M$ be an
indecomposable holonomic $A_n$-module saturated with respect to
$f_1\cdots f_s$. The first step of Algorithm \ref{alg3} also comes
from the fact that the  \v{C}ech complex
${\check{C}}^{\bullet}(f_1,\dots,f_s;M)$ is quasi-isomorphic to
$$
\hskip 5mm 0 \longrightarrow 0 \stackrel{d_0}\longrightarrow
M_{f_1}/M \stackrel{d_1}\longrightarrow \bigoplus_{i=2}^{s}
M_{f_1f_i}/ M_{f_i} \stackrel{d_2}\longrightarrow \cdots
\longrightarrow M_{f_1\cdots f_s}/ M_{f_2\cdots
f_{s}}\longrightarrow 0.
$$

It would be interesting to continue the pruning process for the
whole complex (not just for the components of the characteristic
cycle) in order to find a complete description of a minimal complex
quasi-isomorphic to the \v{C}ech complex.

\vskip 2mm

A canonical \v{C}ech complex is introduced in \cite{MS05} for the case of $I$ being a monomial ideal and $M=R$.
This complex is associated to a minimal free resolution of $R/I$ in the same way as the usual \v{C}ech complex
is associated to the Taylor resolution of $R/I$. The difference with the pruned \v{C}ech complex we propose is that,
roughly speaking, they only prune when the localization map $R_{f_\alpha}\lra R_{f_{\alpha+\varepsilon_i}}$ is the identity.

\end{remark}

\section{Examples}

In this section we want to present some examples where we will apply our algorithm to
compute the characteristic cycle of local cohomology modules. First we have to study localizations $R_f$ of the polynomial ring $R=\bQ[x_1,\dots,x_n]$ at a polynomial $f\in R$. To compute its
characteristic cycle directly one needs to:

\vskip 2mm

\hskip 1cm $\cdot$ Construct a presentation of the $A_n$-module $R_f$,

\hskip 1cm $\cdot$ Compute the characteristic ideal $J(R_f)$,

\hskip 1cm $\cdot$ Compute the primary decomposition of $J(R_f)$ and its corresponding
 multiplicities.

\vskip 2mm

The first two steps require expensive computations in the Weyl
algebra $A_n$ since we have to compute the Bernstein-Sato polynomial
of $f$. For some short examples we can do the job just using the
{\tt Macaulay2} commands {\tt Dlocalize} and {\tt charIdeal}.

\vskip 2mm

Following the approach of this work we have developed some scripts written in {\tt Macaulay 2} that compute and print out the
list of components and the corresponding multiplicities showing up in the characteristic cycles of the localizations $R_f$ in
the examples we present in this section. In fact we develop two different strategies that we may use depending on the examples
we want to treat.

\vskip 2mm

{\it $\cdot$ Single localization:} \hskip 2mm Since the characteristic cycle of $R$ is $CC(R)= T_X^*X$, the characteristic
cycle of $R_f$ is $CC(R_f)=T_X^*X+\Gamma$, where $\Gamma$ is computed according to Theorem \ref{propBMM} so we may compute it
in one step. Notice that the defining ideal of $\Gamma$ may be quite large so computing its primary decomposition can be
expensive.

\vskip 2mm

{\it $\cdot$ Iterative localization:} \hskip 2mm We can apply
Theorem \ref{propBMM} in an iterative way on the components of the
polynomial $f$. This strategy is useful to treat large examples,
since, usually, it leads to computing primary decompositions of
ideals of lower degrees compared to the former strategy.

\vskip 2mm

For the examples we present in this work both strategies can be applied.

\subsection{Local cohomology modules}
\label{subsecLocalCohomology}

Consider the ideal $I\subset R=\bQ[x_1,...,x_6]$ generated by the minors of the matrix
$$
\left(
\begin{array}{lll}
x_1 & x_2 & x_3 \\
x_4 & x_5 & x_6
\end{array}\right).$$

It is a nontrivial problem to show that the local cohomology module $H^3_I(R)$ is nonzero (see \cite[Remark 3.13]{HL90},
\cite{KL02}). For example, {\tt Macaulay 2} runs out of memory before computing this module with the command {\tt localCohom}.
U.~Walther \cite[Example 6.1]{Wa99} gives a complete description of this module using a tailor-made implementation of his
algorithm which is based on the construction of the \v{C}ech complex. The difference with the implementation of the {\tt
Macaulay 2} command is that he uses iterative localization to reduce the complexity in the computation of Bernstein-Sato
polynomials.

\vskip 2mm

Our method makes it possible to prove algorithmically that $H^3_I(R)\neq 0$ from the computation of the characteristic cycles
of the localization modules in the \v{C}ech complex which for this particular example looks like
%\begin{equation} \label{eqCech}
$$(\star) \hskip 5mm 0 \rightarrow R \rightarrow R_{f_1}\oplus R_{f_2}\oplus R_{f_3} \rightarrow R_{f_1f_2}\oplus R_{f_1f_3}\oplus R_{f_2f_3}
\rightarrow R_{f_1f_2f_3} \rightarrow 0,$$
%\end{equation}
where $f_1= x_1x_5-x_2x_4, \hskip 2mm f_2= x_1x_6-x_3x_4 \hskip 2mm
{\rm and} \hskip 2mm f_3= x_2x_6-x_3x_5$.

\vskip 2mm

\begin{remark}
By flat base change we can also deduce the non-vanishing of the local cohomology module $H^3_I(R)$ where $R=k[x_1,...,x_6]$ is
the polynomial ring over any field $k$ of characteristic zero.
\end{remark}

The list of components and their corresponding multiplicities
showing up in the characteristic cycles of the chains in the
\v{C}ech complex $(\star)$ and different from the whole space $X$ that
we get with our script contains 14 elements. A sample entry
is as follows:

\begin{Macaulay2}
\small
\begin{verbatim}
Component = V(ideal (x x  - x x , x x  - x x , x x  - x x ))
                      3 5    2 6   3 4    1 6   2 4    1 5
entries-> HashTable{{0, 1, 2} => 2}
                    {0, 1} => 1
                    {0, 2} => 1
                    {0} => 0
                    {1, 2} => 1
                    {1} => 0
                    {2} => 0
\end{verbatim}
\end{Macaulay2}

\noindent Namely, the component corresponding to the ideal $I$ is
present with multiplicity one in $R_{f_1f_2}$, $R_{f_2f_3}$,
$R_{f_1f_3}$ and with multiplicity two in $R_{f_1f_2f_3}$. The
following is the complete list of 14 components:

\vskip 2mm

$\begin{array}{ccc}
\hskip -1.3cm \hskip -1.8cm A_1 = V(f_1),& \hskip -1.8cm A_2 = V(f_2),& \hskip -1.8cm A_3 = V(f_3),\\
\hskip -1.3cm \hskip -1.1cm B_1 = V(x_3,x_6),& \hskip -1.1cm B_2 = V(x_2,x_5),&\hskip -1.1cm B_3 = V(x_1,x_4),\\
\hskip -1.3cm \hskip -0.6cm C_1 = V(x_3,x_6,f_1),& \hskip -0.6cm C_2 = V(x_2,x_5,f_2),&
\hskip -0.6cm C_3 = V(x_1,x_4,f_3),\\
\hskip -1.3cm D_1 = V(x_1,x_2,x_4,x_5),& D_2 = V(x_1,x_3,x_4,x_6),& D_3 = V(x_2,x_3,x_5,x_6),\\
\hskip -1.3cm \hskip 1cm E = V(x_1,x_2,x_3,x_4,x_5,x_6),& \hskip -2cm F = V(I).&
\end{array}$

\vskip 2mm

\noindent Piecing the results of our computation together we can draw the $3$-hypercube in Figure $1$.

\vskip 2mm

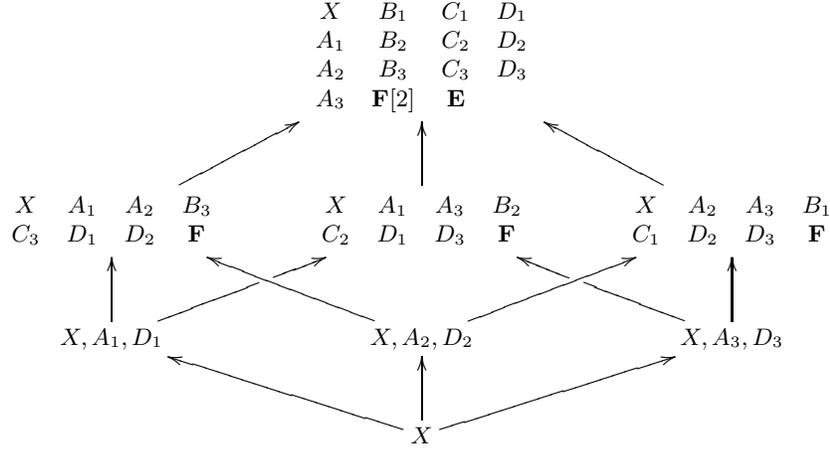
\begin{figure}[h]\label{cocellular2}
{\small $${\xymatrix { &{\begin{array}{cccc}
  X & B_1&C_1&D_1 \\
  A_1&B_2&C_2&D_2\\
  A_2&B_3&C_3&D_3 \\
  A_3&\textbf{F}[2]&\textbf{E}&
\end{array} } &
\\ {\begin{array}{cccc}
  X & A_1 & A_2 &B_3 \\
  C_3&D_1 &D_2 &\textbf{F}
\end{array}} \ar[ur]  & {\begin{array}{cccc}
  X & A_1 & A_3 &B_2 \\
  C_2&D_1 &D_3 &\textbf{F}
\end{array}} \ar[u] & {\begin{array}{cccc}
  X & A_2 & A_3 &B_1 \\
  C_1&D_2 &D_3 &\textbf{F}
\end{array}}  \ar[ul]
\\ X,A_1,D_1 \ar[ur]|\hole \ar[u]&X,A_2,D_2 \ar[ul] \ar[ur]& X,A_3,D_3 \ar[ul]|\hole \ar[u]
\\& X \ar[ul] \ar[u] \ar[ur]& }}$$}

\caption{Components of characteristic cycles for the \^Cech
complex $(\star)$ (multiplicity~$>1$ is specified in square
brackets).}
\end{figure}

To compute the characteristic cycle of the cohomology modules we have to apply Theorem
\ref{T} that has been implemented in the routine {\tt PruneCechComplexCC}.
According to the output
\begin{Macaulay2}
\small
\begin{verbatim}
{} => {}
{0} => {}
{1} => {}
{2} => {}
{0, 1} => {ideal (x x  - x x , x x  - x x , x x  - x x ) => 1}
                   3 5    2 6   3 4    1 6   2 4    1 5
{0, 2} => {}
{1, 2} => {}
{0, 1, 2} => {ideal (x , x , x , x , x , x ) => 1}
                      6   5   4   3   2   1
\end{verbatim}
\end{Macaulay2}
\noindent we get $CC(H^2_I(R)) = T^{\ast}_F X$ and $CC(H^3_I(R)) = T^{\ast}_E
X$. Finally, it is worth to point out that the obtained result is
coherent with the fact that the local cohomology module $H^3_I(R)$
is isomorphic to the injective hull of the residue field
$E_R(R/(x_1,...,x_6))$.

\subsection{Lyubeznik numbers}
Let $R=k[x_1,...,x_n]$ be the polynomial ring over a field $k$ of
characteristic zero. Let $I\subseteq R$ be an ideal and
$\fM=(x_1,...,x_n)$ be the homogeneous maximal ideal. G.~Lyubeznik
\cite{Ly93} has defined a new set of numerical invariants of the
quotient ring $R/I$ by means of the Bass numbers
$$\la_{p,i}(R/I):=\mu_p(\fM,H_I^{n-i}(R)) := {\dm}_k \hskip 1mm
{\Ext}_R^p(k,H_I^{n-i}(R)).$$ These invariants can be described as
the multiplicities of the characteristic cycle of the local
cohomology modules $H_{\fM }^p(H_{I}^{n-i}(R))$ (see \cite{Al02}).
Namely, $$CC(H_{\fM }^p(H_{I}^{n-i}(R)))= \la_{p,i} \hskip 1mm
T^{\ast}_E X $$

\vskip 2mm

Lyubeznik numbers carry interesting topological information of the quotient ring $R/I$ as it is pointed in \cite{Ly93} and
\cite{GS98}. To compute them for a given ideal $I\subseteq R$ and arbitrary $i,p$ we refer to U.~Walther's
algorithm \cite[Algorithm 5.3] {Wa99} even though it has not been implemented yet. When $I$ is a squarefree monomial ideal, a description of
these invariants is given in \cite{Al00}. Some other particular computations may also be found in \cite{GS98} and \cite {Wa01}.

\vskip 2mm

Let $I\subset R=\bQ[x_1,...,x_6]$ be the ideal generated by the minors of the matrix
$$
\left(
\begin{array}{lll}
x_1 & x_2 & x_3 \\
x_4 & x_5 & x_6
\end{array}\right)$$ considered above, i.e.
$I= (x_1x_5-x_2x_4, x_1x_6-x_3x_4, x_2x_6-x_3x_5)$. We want to compute the characteristic cycle of the local cohomology modules
$H_{\fM }^p(H_{I}^{i}(R))$ for $i=2,3$ and $\forall p$ so we have to construct the \v{C}ech complex
%\begin{equation} %\label{eqCech2}
$$(\star \star) \hskip 5mm 0 \rightarrow M \rightarrow \bigoplus_{i=1}^{6} M_{x_i} \rightarrow \cdots  \rightarrow M_{x_1\cdots x_6} \rightarrow 0,$$
%\end{equation}\star \star
where $M$ is either $H_{I}^{2}(R)$ or $H_{I}^{3}(R)$. Then we have to compute the characteristic cycles of the localization
modules and use Theorem \ref{T}.

\vskip 2mm

$\bullet$ For $M=H^3_I(R)$ we know that its characteristic cycle is $T^{\ast}_E X$ so, applying Theorem \ref{propBMM}, the
\v{C}ech complex $(\star \star)$ reduces to the first term. Then,
$$CC(H^0_{\fM }(H^3_I(R))) = T^{\ast}_E X$$ and the other local cohomology modules vanish.

\vskip 2mm

$\bullet$ For $M=H^2_I(R)$ we obtain $$CC(H^2_{\fM }(H^2_I(R))) = T^{\ast}_E X$$ $$CC(H^4_{\fM }(H^2_I(R))) = T^{\ast}_E X$$
and the other local cohomology modules vanish. We are not going to present the complete output with all the components as in Figure $1$ for this case but at least we are going to show the multiplicities of the component $T^{\ast}_E X$ appearing in the
\v{C}ech complex $(\star \star)$ in Figure $2$. We point out that the components $T^{\ast}_E X$ that survive to Algorithm \ref{alg3} belong to $CC(M_{x_1x_2})$ and $CC(M_{x_1x_4x_5x_6})$.

\begin{figure}[h]
{\small $$\hskip -1cm {\xymatrix { &{\begin{array}{cccccccccccccccccccc}
M: \emptyset &M_{x_1}: \emptyset &M_{x_1x_2}:  E& M_{x_1x_2x_3}:  E[2]& M_{x_1x_2x_3x_4}:  E[3]          &M_{x_1x_2x_3x_4x_5}:  E[3]& M_{x_1x_2x_3x_4x_5x_6}:  E[3]\\
&M_{x_2}: \emptyset&M_{x_1x_3}:  E&M_{x_1x_2x_4}:  E&M_{x_1x_2x_3x_5}:  E[3]&M_{x_1x_2x_3x_4x_6}:  E[3]& \\
&M_{x_3}: \emptyset&M_{x_1x_4}:  \emptyset&M_{x_1x_2x_5}:  E&M_{x_1x_2x_3x_6}:  E[3]&M_{x_1x_2x_3x_5x_6}:  E[3]& \\
&M_{x_4}: \emptyset&M_{x_1x_5}:  E&M_{x_1x_2x_6}:  E[3]&M_{x_1x_2x_4x_5}:  E&M_{x_1x_2x_4x_5x_6}:  E[3]& \\
&M_{x_5}: \emptyset&M_{x_1x_6}:  E&M_{x_1x_3x_4}:  E&M_{x_1x_2x_4x_6}:  E[3]&M_{x_1x_3x_4x_5x_6}:  E[3]& \\
&M_{x_6}: \emptyset&M_{x_2x_3}:  E&M_{x_1x_3x_5}:  E[3]&M_{x_1x_2x_5x_6}:  E[3]&M_{x_2x_3x_4x_5x_6}:  E[3]& \\
&&M_{x_2x_4}:  E&M_{x_1x_3x_6}:  E&M_{x_1x_3x_4x_5}:  E[3]&& \\
&&M_{x_2x_5}:  \emptyset &M_{x_1x_4x_5}:  E&M_{x_1x_3x_4x_6}:  E&& \\
&&M_{x_2x_6}:  E&M_{x_1x_4x_6}:  E&M_{x_1x_3x_5x_6}:  E[3]&& \\
&&M_{x_3x_4}:  E&M_{x_1x_5x_6}:  E[3]&M_{x_1x_4x_5x_6}:  E[3]&& \\
&&M_{x_3x_5}:  E&M_{x_2x_3x_4}:  E[3]&M_{x_2x_3x_4x_5}:  E[3]&& \\
&&M_{x_3x_6}:  \emptyset &M_{x_2x_3x_5}:  E&M_{x_2x_3x_4x_6}:  E[3]&& \\
&&M_{x_4x_5}:  E&M_{x_2x_3x_6}:  E&M_{x_2x_3x_5x_6}:  E&& \\
&&M_{x_4x_6}:  E&M_{x_2x_4x_5}:  E&M_{x_2x_4x_5x_6}:  E[3]&& \\
&&M_{x_5x_6}:  E&M_{x_2x_4x_6}:  E[3]&M_{x_3x_4x_5x_6}:  E[3]&& \\
&&&M_{x_2x_5x_6}:  E&&& \\
&&&M_{x_3x_4x_5}:  E[3]&&& \\
&&&M_{x_3x_4x_6}:  E&&& \\
&&&M_{x_3x_5x_6}:  E&&& \\
&&&M_{x_4x_5x_6}:  E[2]&&&
                                   \end{array}
}}}$$}

\caption{Component $T^{\ast}_E X$ appearing in the \v{C}ech complex $(\star \star)$ (multiplicity $>1$ is specified in square
brackets).}

\end{figure}

\vskip 2mm

Using the properties that Lyubeznik numbers satisfy (see \cite[Section 4]{Ly93}), we can collect the multiplicities in a
triangular matrix as follows:
$$\Lambda(R/I)=\begin{pmatrix}
  0 & 0 & 0 & 1 & 0 \\
   & 0 & 0& 0 & 0 \\
   &  & 0 & 0 & 1 \\
   &  &  & 0 & 0 \\
   &  &  &  & 1
\end{pmatrix}$$
The complex variety $V$ defined by $I$ has an isolated singularity at the origin. The singular cohomology groups of $V$ with complex coefficients and support at the origin can be described
from Lyubeznik numbers (see \cite{GS98}). In our case we get

\vskip 2mm

\hskip 1cm $1=\lambda_{4,4}= \dim_{\bC} H^8_{\{0\}}(V,\bC)$,

\hskip 1cm $1=\lambda_{2,4}= \dim_{\bC} H^6_{\{0\}}(V,\bC)$,

\hskip 1cm $1=\lambda_{0,3}= \dim_{\bC} H^3_{\{0\}}(V,\bC)$.

\section{Conclusion and possible developments}

We have shown that characteristic cycles of local cohomology modules
can be be computed by algorithm operating in commutative polynomial
rings as opposed to the direct computation of these modules, which
requires Gr\"{o}bner bases technique in noncommutative Weyl
algebras.

\vskip 2mm

The computational engine of our method is primary decomposition,
from which we extract only geometrical information. This prompts a
natural interest in \emph{numerical primary decomposition}, an
algorithm that would produce just that -- the descriptions of
reduced components and their multiplicities -- by means of numerical
computations.

\end{document}